\UseRawInputEncoding 

\documentclass[12pt,eqno]{article}

\usepackage{amssymb, amsfonts, eucal, amsmath}

\usepackage{amsfonts,amssymb,amsmath}

\newtheorem{lem}{Lemma}[section]

\def\ol {\overline }

\def\al {\alpha }
\def\si {\sigma }

\newtheorem{theorem}{Theorem}

\newtheorem{remark}{Remark}

\def\proof{\medskip\noindent{\it Proof.} }

\newcommand{\sect}[1]{\section{#1} \setcounter{equation}{0} }

\begin{document}


\begin{center}

{\Large\bf Degree of nearly comonotone approximation \\
\vskip 0.2cm
of periodic functions }\label{contr comono otr}

\vskip 0.2cm

{ German Dzyubenko}

\vspace{0.2cm}

\footnotesize{\it
Institute of Mathematics NAS of Ukraine, \\
01024 Ukraine, Kyiv-4, 3, Tereschenkivska st.,\\
\rm dzyuben@gmail.com
}

\setcounter{footnote}{0}
\renewcommand{\thefootnote}{ }


\vskip 0.2 cm

\begin{minipage}{14cm}

{\footnotesize\bf Abstract.} {\footnotesize
Let a $2\pi$-periodic function $f\in\Bbb C$ changes its monotonicity at a finitely even number of points $y_i$ of the period.
The
degree of approximation of this $f$ by trigonometric polynomials which are comonotone
with it, i.e. that change their monotonicity exactly
at the points $y_i$ where $f$ does, is restricted by $\omega_2(f,\pi/n)$ (with a
constant depending on the location of these $y_i$). Recently, we proved that
relaxing the comonotonicity requirement in intervals of length proportional to $\pi/n$ about the
 points $y_i$ (so called nearly comonotone approximation) allows the polynomials to achieve the approximation rate of $\omega_3$.
By constructing a counterexample, we show
here that even with the relaxation of the requirement of comonotonicity for the polynomials on sets
with measures approaching $0$ (no matter how slowly or how fast) $\omega_4$ is not reachable.
}

\end{minipage}

\end{center}

\sect{Introduction}

$\bf 1^\circ .$ Let $\Bbb C:=\Bbb C_{\Bbb R}$ and $\Bbb C^r$ denote, respectively,
 the space of continuous $2\pi$-periodic
functions $f:\Bbb R\rightarrow\Bbb R$, and that of $r$-times continuously differentiable functions
on the real axis $\Bbb R,$
equipped
with the uniform norm
$$
\|f\|:=\| f\|_{\Bbb R}:=\max_{x\in \Bbb R}|f(x)|;
$$
$\Bbb T_n,\ n\in\Bbb N,$ is the space of trigonometric polynomials
 $P_n(x)=a_0+\sum_{j=1}^{n}(a_j\cos
jx+b_j\sin jx) $ of degree $\le n$ (of order $\le 2n+1 $) with $a_j,b_j\in\Bbb R.$ By
$$
E_n(f):=\inf\limits_{P_n\in\Bbb T_n}\left\|f-P_n\right\|
$$
we denote
the value (error) of the best uniform approximation of the function $f$ by polynomials
$P_n\in\Bbb T_n.$
For any bounded  $2\pi$-periodic function $f,$ and
$k\in N,$ the $k$-th symmetric difference of $f$ at the point
$x$ with the step $h\ge 0$ is defined as
$$
\Delta^k_hf(x):=\sum_{i=0}^k(-1)^{k-i}\binom kif(x-\frac k2h+ih),\quad x\in[-\pi,\pi],
$$
and the (ordinary) $k$-th (order) modulus of continuity (or smoothness) of  $f\in\Bbb C$ is defined as
$$
\omega_k(f,t):=\sup_{0<h\le t}\sup_{x\in[-\pi,\pi]}|\Delta^k_hf(x)|.
$$

We recall the classical Jackson-Zygmund-Akhiezer-Stechkin estimate
(obtained by Jackson for $k=1$ \cite{jack,jack2}, Zigmund for $k=2,\ \omega_2(f,t)\le t,$ \cite{zygm}, Akhiezer for $k=2$ \cite{ahie} and Stechkin for $k\ge
3$ \cite{stec}):
{\it If a function $f\in\Bbb C,$ then
\begin{equation}\label{classical}
E_n(f) \le c(k)\, \omega_k \left(f,\pi/{n}\right),\quad
n\in\Bbb N,
\end{equation}
where $c(k)$ is a constant that depends only on $k,$
} for details see, for example, \cite[Section 4]{dzya}. And hence, in particular, {\it if $f\in\Bbb C^r,\ r\in\Bbb N,$ then
\begin{equation}\label{classical1}
E_n(f) \le \frac{c(r,k)}{n^r}\, \omega_k \left(f^{(r)},\pi/{n}\right),\quad
n\in\Bbb N.
\end{equation}
}

In 1968 Lorentz and Zeller \cite{LZ, LZ1969} proved a bell-shaped analogue of the estimate \eqref{classical} with $k=1$, that is,
for approximation of bell-shaped (i.e. even and nonincreasing on $[0,\pi]$) functions from $\Bbb C$ by bell-shaped
polynomials from $\Bbb T_n$, and thus gave rise to the search for its other analogues, i.e. with other restrictions on the shape of the function and polynomials such as piecewise positivity, monotonicity, convexity
 (now this is called Shape Preserving Approximation, see, for example, the surveys of  Kopotun, Leviatan, Prymak, Shevchuk \cite{klps,kls}).

\vskip 0.2cm

\noindent
$\bf 2^\circ .$
In the paper of Pleshakov \cite{pl}, and in \cite{dzpl} a comonotone analogue of  \eqref{classical} is proved with $k=1$ and $k=2,$ respectively. To write it we need some notations.
Let on  $ [-\pi,\pi ) $ there are $ 2s,\ s\in\Bbb N, $  fixed points
$$ y_i:\ -\pi\le y_{2s}<y_{2s-1}<...<y_{1}<\pi , $$
and for the rest $ i\in\Bbb
Z,  $ the points $ y_i $ are defined by the equality $ y_i=y_{i+2s}+2\pi $
(i.e. $y_0=y_{2s}+2\pi ,..., y_{2s+1}=y_1-2\pi ,... $).
Let $ Y_{s}:=\{
y_i\}_{i\in\Bbb Z}.$ We denote by $\Delta^{(1)}(Y_{s}) $ the collection of all functions $
f\in\Bbb C $ that are nondecreasing  on $[y_1,y_0],$ nonincreasing  on $[y_2,y_1],$
nondecreasing on $[y_3,y_2]$ and so on. Thus, if $f\in\Bbb C^1$ then $f\in\Delta^{(1)}(Y_{s})\, \Leftrightarrow \,
f'(x)\Pi(x)\ge 0,\ x\in\Bbb R,$ where
$$
\Pi(x):=\Pi(x,Y_{s}):=\prod_{i=1}^{2s}\sin\frac{1}{2}(x-y_i)
$$
$(\Pi(x)>0,\ x\in(y_1,y_0),\ \Pi\in\Bbb T_s).$
The functions from  $\Delta^{(1)}(Y_{s})$ are  called  piecewise monotone or {\it comonotone} (each other or between themselves), and the approximation of them by polynomials also from $\Delta^{(1)}(Y_{s})$ is called comonotone approximation.
Let
$$
E_n^{(1)}(f,Y_{s}):=\inf\limits_{P_n\in\Bbb
T_n\cap\Delta^{(1)}(Y_{s})}\left\|f-P_n\right\|
$$
is the value (error) of the best uniform approximation of the function $f$ by polynomials
$P_n\in\Bbb T_n\cap\Delta^{(1)}(Y_{s}).$

So, it is proved in \cite{pl} and \cite{dzpl} that:
 {\it  If a function $f\in\Delta^{(1)}(Y_{s}), $ then there exists a constant $ N(Y_{s}),$ which depends only on $\min\limits_{i=1,...,2s}\,\{y_{i}-y_{i+1}\},$ such that
\begin{equation}\label{om2}
E_n^{(1)}(f,Y_{s}) \le c(s)\,
\omega_k (f,{\pi}/{n}), \quad k=1,2,\quad n\ge N(Y_{s}),
\end{equation}
where  $c(s):=\max\limits_{k=1,2}c(s,k)$ is a constant depending only on $s.$
}

Note that the estimate \eqref{om2'} below  is a simple consequence of \eqref{om2} and the Whitney inequality
 \cite{whit} $\Vert f-f(0)\Vert\le c\,\omega_k(f,2\pi)$ with an absolute constant $c$ and $k\in\Bbb N,$ (since $f(0)\in\Delta^{(1)}(Y_{s})\cap\Bbb T_0$ and at least twice interpolates $f$ on the period):
{\it
 If $f\in\Delta^{(1)}(Y_{s}),$ then
\begin{equation}\label{om2'}
E_n^{(1)}(f,Y_{s}) \le C(Y_{s})\, \omega_2 \left(f,\pi/{n}\right),\quad n\in\Bbb N,
\end{equation}
 where $C(Y_{s})$ is a constant depending only on $\min\limits_{i=1,...,2s}\,\{y_{i}-y_{i+1}\}$.
 }

Moreover, Pleshakov and Tyshkevich  \cite{pl_dis, plty} using considerations of the papers of
Shvedov \cite {shve1} and DeVore, Leviatan, Shevchuk \cite {dls}, {\it for each $n\in\Bbb N,$ constructed
a function $g_n(x)=g_n(x,Y_{s})\in\Delta^{(1)}(Y_{s})$ such that}
$$
E_n^{(1)}(g_n,Y_{s})\ge C(Y_{s})\,\sqrt{n}\, \omega_3 \left(g_n,\pi/n\right)
$$
for some $C(Y_{s})>0.$
In other words,
for each $n\in\Bbb N,$ they found a function from $\Delta^{(1)}(Y_{2}),$ for which the inequality
\eqref{om2'} is invalid with $\omega_k, \ k>2,$ i.e. it is impossible to improve it in the order of the modulus of smoothness (unlike \eqref {classical} and \eqref {classical1} of approximation without restrictions, which hold for all $ k\in\Bbb N$).
In \cite{dz_con}, one such function is constructed (for all $n\in\Bbb N$), i.e. it is proved that
{\it in the set $\Delta^{(1)}(Y_{s})$ there exists a function $g(x)=g(x,Y_{s})$ such that }
$$
\limsup_{n\rightarrow\infty}\frac{E_n^{(1)}(g,Y_{s})}{\omega_3(g,\pi/n)}=\infty.
$$

\vskip 0.1cm

\noindent
$\bf 3^\circ .$
However, from the result of Leviatan and Shevchuk \cite{lesh1} on approximation on a segment by algebraic polynomials, we knew that relaxing the comonotonicity requirement in small intervals about the points $y_i$
(so called {\it nearly} comonotone approximation) allows the polynomials to achieve one additional approximation rate, and it was proved in \cite{dznea} that:
{\it If
 $f\in\Delta^{(1)}(Y_{s})$ then there exists a constant  $N(Y_{s}),$ which depends only on
 $\min\limits_{i=1,...,2s}\,\{y_{i}-y_{i+1}\},$ and for every $n\ge N(Y_{s})$ there exists a polynomial
$P_n\in\Bbb T_{cn}$ such that
$$
P_n'(x)\Pi(x)\ge 0,\quad x\in\Bbb R\setminus \cup_{i\in\Bbb Z}\left(y_i-\pi/n,y_i+\pi/n\right),
$$
\begin{equation}\label{estimate}
E_n^{(1)}(f,Y_{s})\le
\|f-P_n\|\le c(s)\,\omega_3(f,\pi/n),
\end{equation}
where  $c$ and $ c(s)$ are constants depending only on   $s.$
}

Like \eqref{om2'} the estimate \eqref{estimate2} below  follows from \eqref{estimate}  and the Whitney inequality
 \cite{whit}:
{\it If
 $f\in\Delta^{(1)}(Y_{s})$ then for every $n\in\Bbb N$ there exists
$P_n\in\Bbb T_n$ such that
$$
P_n'(x)\Pi(x)\ge 0,\quad x\in\Bbb R\setminus \cup_{i\in\Bbb Z}\left(y_i-c/n,y_i+c/n\right),
$$
\begin{equation}\label{estimate2}
E_n^{(1)}(f,Y_{s})\le
\|f-P_n\|\le C(Y_{s})\,\omega_3(f,\pi/n),
\end{equation}
where  $c$
and $C(Y_{s})$
are constants depending only on $s$ and
$\min\limits_{i=1,...,2s}\,\{y_{i}-y_{i+1}\}.$
}

\vskip 0.1cm

\noindent
$\bf 4^\circ .$
In this paper we prove Theorem \ref{thm1} below. For any $\epsilon >0$ and $f\in\Delta^{(1)}(Y_{s})$ let
$$
E_n^{(1)}(f,\epsilon ,Y_{s}):=\inf_{P_n}\|f-P_n\|,
$$
where the infimum is taken over all polynomials $P_n\in\Bbb T_n$ satisfying
$$
\text{meas} \bigl(\{x:\, P_{n}'(x)\Pi(x,Y_{s})\ge 0\}\cap [-\pi,\pi]\bigl)\ge 2\pi -\epsilon .
$$
(Obviously, $E_n^{(1)}(f,Y_{s})=E_n^{(1)}(f,0,Y_{s})$.)

\begin{theorem}\label{thm1}
 For each $s\in\Bbb N$ and each sequence $\bar\epsilon=\{\epsilon_n\}_{n=1}^{\infty}$ of nonnegative numbers tending to $0,$ there exist a set $Y_{s}$ and a function $f=f_{\bar\epsilon}\in\Delta^{(1)}(Y_{s})$ such that
\begin{equation}\label{main}
\limsup_{n\rightarrow\infty}\frac{E_n^{(1)}(f,\epsilon_n,Y_{s})}{\omega_4(f,\pi/n)}=\infty.
\end{equation}
\end{theorem}

\begin{remark}\label{rem1}
The considerations in the proof of Theorem \ref{thm1} are inspired by their algebraic analogue, namely, by the paper of  Leviatan and Shevchuk \cite{lesh2}, and their joint paper with DeVore \cite {dls}. Partially these considerations were also used in \cite[Example 3.2]{dgs}, \cite{dz_con} and \cite[Theorem 4.1]{dvy}.
\end{remark}
\begin{remark}\label{rem2}
In the algebraic case, i.e. in \cite{lesh2}, a stronger result is obtained, namely, the algebraic analogue of Theorem \ref{thm1} is stated there for an arbitrary $Y_{s}\subset (-1,1),$ and this is natural, because no need to regard the periodicity.  It seems that for an arbitrary $Y_{s}$ in the trigonometric case, it is necessary to construct a different counterexample than in the proof of Theorem \ref{thm1}.
\end{remark}


\sect{  Construction of a counterexample \\ (Proof of Theorem \ref{thm1})}

As we said, this paper is inspired by the work of Leviatan and Shevchuk \cite {lesh2} and their joint work with DeVore \cite {dls}, so now we will follow along the lines of \cite {lesh2} and \cite {dls} with the following main differences: instead of the Chebyshev and truncated Chebyshev (algebraic) polynomials we will use (have to use) others (trigonometric), to preserve the periodicity we will use not an arbitrary $Y_{s}$ rather a specific $Y^*_{s}$ (with equidistant $y_i$'s)
and we will change the definition of the function $g_{j,M}$ in \cite {lesh2} for a more precise description of the continuity of the function $g$ (see \cite [after (2.7)] {lesh2}).

\vskip 0.2cm

\noindent
$\bf 1^\circ .$
Let $Y^*:=Y^*_{s}=\{y_i:= \pi+\frac{\pi}{2s}-i\frac{\pi}{s},\ i=1,...,2s\},\ s\in\Bbb N$ is even for definiteness, $[y_{s+1},y_{s}]$ is the interval containing $0$ and $I_0:=[-\frac{\pi}{3s},\frac{\pi}{3s}]$ is its middle two thirds (central four parts of six), so that, in particular, $\Pi(0)>0,\ \|\Pi\|=\|\Pi\|_{[-\pi,\pi]}=\frac{1}{2^{2s-1}},$ and
\begin{equation}\label{minpi}
\min\limits_{x\in I_0}|\Pi(x)|=\Pi(\pm\frac{\pi}{3s})=\frac{1}{2^{2s}},
\end{equation}
here and in the sequel $\Pi=\Pi(x)=\Pi(x,Y^*)$ (one can verify the both equalities
using math software).

For $\nu > 10s$ and even for definiteness, let
$$
t_{\nu}(x):=\cos\nu\,x,
$$
with the extrema $x_j:=\pi-\frac{j\pi}{\nu}, \ -\infty<j<\infty.$  Given $0<b<\frac12,$ we take two points on both side of $x_j,$ namely, set $x_{j,l}:=\pi-\frac{(j+b)\pi}{\nu}$ and $x_{j,r}:=\pi-\frac{(j-b)\pi}{\nu}.$ Then $$
|t_{\nu}(x_{j,l})|=|t_{\nu}(x_{j,r})|=\cos\pi b,
$$
and
\begin{equation}\label{r-l}
x_{j,r}-x_{j,l}=2\pi\frac{b}{\nu}.
\end{equation}
We also need the truncated $t_{\nu},$ namely,
$$
t^*_{\nu}(x):=t^*_{\nu,b}(x):=\left\{ \aligned
-\cos\pi b,\ \ & \text{if} \ \ t_{\nu}(x)<-\cos\pi b,\\
t_{\nu}(x),\ \ & \text{otherwise}.
\endaligned
\right.
$$
Put
$$
t_{\nu,b}:={2^{2s}}(t_{\nu}+\cos\pi b)\Pi,\ \ \text{and}\ \
\ol{t}_{\nu,b}:={2^{2s}}(t^*_{\nu,b}+\cos\pi b)\Pi.
$$
It is readily seen that
\begin{align}\label{tnu}
\|\ol{t}_{\nu,b}\|\le\|{t}_{\nu,b}\|&=\|{t}_{\nu,b}\|_{[-\pi,\pi]}\nonumber \\
&\le 2\cdot {2^{2s}}\|\Pi\|_{[-\pi,\pi]}=4.
\end{align}

Now, we define a function $g=g_{\nu,b}$ on $[-\pi,\pi]$ by defining it on each interval $[x_{2j+1},x_{2j-1}]\cap[-\pi,\pi],\ j=0,...,\nu,$ (containing one positive extrema of $t_{\nu}$ at its center and at an endpoint for $j=0,\nu$). Note that there are at most $4s$ intervals $[x_{2j+1},x_{2j-1}]$ such that $Y^*\cap(x_{2j+2},x_{2j-2})\ne\emptyset$ (roughly speaking either on both sides of each $y_i$ or directly on it) we call them type I and put
$$
g(x):=\ol{t}_{\nu,b}(x),\quad x\in [x_{2j+1},x_{2j-1}].
$$
All other intervals $[x_{2j+1},x_{2j-1}]\cap[-\pi,\pi],\ j=0,...,\nu,$ are of type II and satisfy $Y^*\cap[x_{2j+1},x_{2j-1}]=\emptyset.$ For these intervals the definition of $g$ requires some preliminaries. Note that for each $1\le i\le 2s$ and $x\in [x_{2j+1},x_{2j-1}],$ in view of $\sin t>\frac3{\pi}t,$ for $0<t<\frac{\pi}{6},$ we have
\begin{equation}\label{231}
\frac{1}{\pi}\le\frac{\sin\frac{x-y_i}{2}}{\sin\frac{x_{2j}-y_i}{2}}\le\frac{2\pi}{3},
\end{equation}
which in turn implies
\begin{equation}\label{pipi}
\left(\frac{1}{\pi}\right)^{2s}\le\frac{\Pi(x)}{\Pi(x_{2j})}\le \left(\frac{2\pi}{3}\right)^{2s}.
\end{equation}
Indeed, without loss of generality assume that $y_i\le x_{2j+2}$. Then if $x>x_{2j},$
then
$$
\frac{2}{\pi}\le\frac{\sin\frac{\pi}{\nu}}{\sin\frac{3\pi}{2\nu}}\le\frac{\sin\frac{x-y_i}{2}}{\sin\frac{x_{2j}-y_i}{2}}
\le\frac{\sin\frac{2\pi}{\nu}}{\sin\frac{\pi}{\nu}} \le \frac{2\pi}{3},
$$
while if $x<x_{2j},$
then
$$
\frac{\pi}{2}\ge\frac{\sin\frac{3\pi}{2\nu}}{\sin\frac{\pi}{\nu}}\ge\frac{\sin\frac{x-y_i}{2}}{\sin\frac{x_{2j}-y_i}{2}}
\ge\frac{\sin\frac{\pi}{2\nu}}{\sin\frac{3\pi}{2\nu}} \ge \frac{1}{\pi}.
$$
Thus \eqref{231} is proved
(arithmetic with math tools gives
$
\frac{2}{\pi}\le\frac{\sin\frac{x-y_i}{2}}{\sin\frac{x_{2j}-y_i}{2}}\le\frac{6}{\pi}
$
).
Now,  in view of \eqref{pipi}, it follows that there is a constant $B < 1/2$,
depending only on $s,$ such that if $b < B,$ then
\begin{equation}\label{piint}
\Pi(x_{2j})\int_{x_{2j+1}}^{x_{2j-1}}t_{\nu,b}(x)dx>0.
\end{equation}
Let $m_j:=\|\ol{t}_{\nu,b}\|_{[x_{2j+1},x_{2j-1}]}.$ For each $0\le M\le 4,$ where $4$ is from \eqref{tnu}, and $x\in [x_{2j+1},x_{2j-1}],$ we set
$$
g_{j,M}(x):=\left\{ \aligned
\frac{M}{m_j}\ol{t}_{\nu,b}(x),\ \ & \text{if} \ \ m_j>M,\\
\ol{t}_{\nu,b}(x),\ \ & \text{otherwise}.
\endaligned
\right.
$$
We now define the function $g$ in intervals $[x_{2j+1},x_{2j-1}]$ of type II where
$\Pi(x_{2j})>0$.  The
case where $\Pi(x_{2j})<0$ requires obvious modifications in its proof and is left without specification. Since
$$
\int_{x_{2j+1}}^{x_{2j-1}}(t_{\nu,b}(x)-g_{j,4}(x))dx=\int_{x_{2j+1}}^{x_{2j-1}}(t_{\nu,b}(x)-\ol{t}_{\nu,b}(x))dx<0,
$$
and by \eqref{piint},
$$
\int_{x_{2j+1}}^{x_{2j-1}}(t_{\nu,b}(x)-g_{j,0}(x))dx=
\int_{x_{2j+1}}^{x_{2j-1}}t_{\nu,b}(x)dx>0,
$$
hence a constant $0<M_j<4$ exists such that
\begin{equation}\label{intint}
\int_{x_{2j+1}}^{x_{2j-1}}t_{\nu,b}(x)dx=\int_{x_{2j+1}}^{x_{2j-1}}g_{j,M_j}(x)dx.
\end{equation}
Thus we define
$$
g(x):=g_{j,M_j}(x),\quad x\in[x_{2j+1},x_{2j-1}].
$$
This completes the definition of $g=g_{\nu,b},$ and
in view of the fact that $g(x_{2j\pm 1})=\ol{t}_{\nu,b}(x_{2j\pm 1})=0,$ for all
$x_{2j}\in[-\pi,\pi],\ j=0,...,\nu,$ it follows that $g$ is continuous in $[-\pi,\pi]$ and hence in $\Bbb R$, and $2\pi$-periodic. We also note that
\begin{equation}\label{gpi}
g(x)\Pi(x)\ge 0,\quad x\in\Bbb R.
\end{equation}
Next, for any interval $[x_{2j+1},x_{2j-1}]$ of either type I or II, from \eqref{r-l} and \eqref{piint} it is follows that
\begin{equation}\label{t-g}
\aligned
\int_{[x_{2j+1},x_{2j-1}]}&|t_{\nu,b}-g(x)|dx \\
&\le 2\left(\int_{x_{2j+1}}^{x_{2j+1,r}}+\int_{x_{2j-1,l}}^{x_{2j-1}}\right) |t_{\nu,b}|dx \\
&\le 2((x_{2j+1,r}-x_{2j+1})+(x_{2j-1}-x_{2j-1,l}))2^{2s}(1-\cos\pi b)\|\Pi\| \\
&<16\pi\frac{b}{\nu}5b^2=:C_0\frac{b^3}{\nu},
\endaligned
\end{equation}
where we used the inequality
$
1-\cos\pi b=2\sin^2\frac{\pi b}{2}<5b^2.
$
Here and in the sequel we denote by $C_0,\ C_1,...,C_9,$
different constants which may
depend only on $s.$

Given $n \ge 1$ and $0<b<B,$ let $\nu=[b^{3/4}n]+13,$ where $[\al]$ denotes
the largest odd integer not exceeding $\al.$ Put
$$
T_{\nu,b}:=\int_0^x t_{\nu,b}(u) du\quad \text{and}\quad f_{\nu,b}:=\int_0^x g_{\nu,b}(u) du.
$$
Since $t_{\nu,b},\ g_{\nu,b},\ \Pi$ and $t_{\nu}$ are even periodic functions and besides $\Pi$ and $t_{\nu}$ have equidistant zeros, then it follows from \eqref{gpi} that
$f_{\nu,b}\in\Delta^{(1)}(Y^*)$. Also, in view of \eqref{t-g}, the following estimate holds
\begin{equation}\label{t-f}
\|T_{\nu,b}-f_{\nu,b}\|=\|T_{\nu,b}-f_{\nu,b}\|_{[-\pi,\pi]}\le C_0(4s+2)\frac{b^3}{\nu} =:C_1\frac{b^3}{\nu}\le C_1\frac{b^{9/4}}{n},
\end{equation}
where we have counted the $4s$ intervals of type I, and possibly two intervals of type II,
namely, those containing $0$ and $x$, respectively.

Set $\ol{x}_{j,l}:=\pi-\frac{(j+b/2)\pi}{\nu}$ and $\ol{x}_{j,r}:=\pi-\frac{(j-b/2)\pi}{\nu},$ and note
\begin{equation}\label{xx}
\ol{x}_{j,r}-x_j=x_j-\ol{x}_{j,l}=
\frac{\pi b}{2\nu},\quad j\in\Bbb Z.
\end{equation}
Let $j$ be odd. Since $\sin b\pi/4> 3b/4$ for $b$ satisfying $b\pi/4<\pi/6,$ then for $x\in[\ol{x}_{j,l},\ol{x}_{j,r}]$, with \eqref{minpi}, the following inequality holds
\begin{equation}\label{t1}
\aligned
T_{\nu,b}'(x)=t_{\nu,b}(x)&< -\cos b\pi/2 +\cos b\pi =-2\sin b\pi/4\sin 3b\pi/4\\
&<-2\frac{3b}{4}\frac{3b}{2}=-\frac{9b^2}{4}.
\endaligned
\end{equation}

It follows by the Bernstein inequality
that
\begin{equation}\label{t4}
\aligned
\|T_{\nu,b}^{(4)}\|&=\|t_{\nu,b}^{(3)}\|_{[-\pi,\pi]} \\
&\le \|t_{\nu,b}\|_{[-\pi,\pi]}\,(\nu+2s)^3\le 4(1+2s)^3\nu^3=:C_2\nu^3.
\endaligned
\end{equation}
Hence, by \eqref{t-f},
\begin{equation}\label{om4}
\aligned
\omega_4(f_{n,b},{\pi}/{n})&\le\omega_4(f_{n,b}-T_{\nu,b},{\pi}/{n})+\omega_4(T_{\nu,b},\pi/{n}) \\
&\le (2\pi)^4\|f_{\nu,b}-T_{\nu,b}\|_{[-\pi,\pi]}+\frac{\pi^4}{n^4}\|T_{\nu,b}^{(4)}\|\\
&\le (2\pi)^4\,C_1\frac{b^{9/4}}{n}+\frac{\pi^4\,C_2\nu^3}{n^4}\le C_3\frac{b^{9/4}}{n},
\endaligned
\end{equation}
where the value of $\nu$ was also used.

{\lem\label{lem0} For any interval $J\subseteq I_0$ there exists an absolute constant $C_4$ such that, if a polynomial $P_n$ satisfies
\begin{equation}\label{p'0}
P_{n}'(x)\ge 0,\quad x\in J\setminus E,
\end{equation}
where $E\subseteq I_0$ is any
 measurable set, then
\begin{equation}\label{fnb-pn}
\|f_{n,b}-P_{n}\|_J\ge \frac{b^2|J|}{\pi n}-\frac{C_4}{n}\left(b^{9/4}+b|E|+\frac{b^{5/4}}{n}\right).
\end{equation}
}

\proof
Let $J_0$ denotes the middle third of $J.$ We consider two cases. First
assume that $J_0$
contains at most one of the $x_j$'s such that $f_{n,b}'(x_j)=0,$ then by the definition of $\nu$
the following inequality holds
$$
|J|< C_5\frac{\pi}{\nu}< C_6\frac{b^{-3/4}}{n},
$$
and hence
\begin{equation}\label{fnb-pn0}
\|f_{n,b}-P_{n}\|_J\ge 0> \frac{b^2|J|}{n}-C_6\frac{b^{5/4}}{n^2}.
\end{equation}
On the other hand, if $J_0$ contains at least two such extrema, then it contains
at least $\frac{4}{\pi}C_7\nu|J|$ of them for some constant $C_7.$
These extrema satisfy \eqref{xx}, and about half of them (and at least one) have odd indices, then
together with \eqref{t1} we conclude that
\begin{equation}\label{meas}
\text{meas} \bigl(J_0\cap\{x:\ T_{\nu,b}'(x)<-\frac{9b^2}{4}\}\bigl)\ge \frac12 \frac{4}{\pi}\frac{\pi b}{2\nu}C_7\nu|J|=C_7b|J|.
\end{equation}
Now, if $C_7b|J|\le|E|,$ then
\begin{equation}\label{fnb-pn1}
\|f_{n,b}-P_{n}\|_J\ge 0\ge \frac{b^2|J|}{n}-\frac{b|E|}{nC_7}.
\end{equation}
Otherwise, by \eqref{meas}, there is a point $x^*\in J_0\setminus E$ such that
$$
T_{\nu,b}'(x^*)<-\frac{9b^2}{4}.
$$
Hence, \eqref{p'0} yields
$$
\frac{9b^2}{4}\le P_n'(x^*)-T_{\nu,b}'(x^*)\le \frac{2\pi}{|J|}n\| P_n-T_{\nu,b}\|_J,
$$
where we used a special case of I.I.Privalov's Theorem (see in \cite[p.96-98]{priv}) with sharp estimates proved in
\cite[Lemma 3.2]{mls}: {\it for every $R_n\in\Bbb T_n$ and each $0<h\le\pi$ the following inequality holds
$$
\|R_n'\|_{[-h/2,h/2]}\le\frac{n}{\sin\frac {h}{2}}\|R_n\|_{[-h,h]}\le\frac{\pi n}{h} \|R_n\|_{[-h,h]}.
$$}
Therefore it follows from \eqref{t-f} that
\begin{equation}\label{pn-t}
\aligned
\frac{b^2|J|}{\pi n}&\le\frac{9b^2}{4n}\frac{|J|}{2\pi}\le\| P_n-T_{\nu,b}\|_J \\
&\le \| P_n-f_{n,b}\|_J+\| f_{n,b}-T_{\nu,b}\|_J\\
&\le \| P_n-f_{n,b}\|_J+C_1\frac{b^{9/4}}{n}.
\endaligned
\end{equation}
Taking $C_4:=\max\{C_6,1/C_7,C_1\},$ \eqref{fnb-pn} follows from combining \eqref{fnb-pn0}, \eqref{fnb-pn1} and \eqref{pn-t}.
Lemma \ref{lem0} is proved.

\vskip 0.2cm

\noindent
$\bf 2^\circ .$
For a given sequence $\bar\epsilon=\{\epsilon_n\}$ we define $f=f_{\bar\epsilon}.$ Let $b_n:=(\max\{\epsilon_n^2,\frac1n\})^{2/5},$ set $d_0:=1,$ and let
$$
d_j:=\frac{b_{n_j}^{9/4}}{n_j}d_{j-1}=\prod_{\mu=1}^{j}\frac{b_{n_{\mu}}^{9/4}}{n_{\mu}},\quad j\ge 1,
$$
where the sequence $\{n_{\mu}\}$
is defined by induction as follows. First, we choose $n_1$ so large
that $b_{n_1}<\min\{B,|I_0|^8\}$
 so that it satisfies \eqref{minbj} below and let $J_0:=I_0.$ Suppose that
$\{n_1,...,n_{\si-1}\}$ and $J_{\si-2}\subseteq J_{\si-3}\subseteq ...\subseteq J_{0},\ \si\ge 2,$
 have been defined. Then put
 $$
 F_{\si-1}:=\sum\limits_{j=1}^{\si-1}d_{j-1}f_{n_j,b_{n_j}},
 $$
and let $J_{\si-1}$ be an interval such that $J_{\si-1}\subseteq J_{\si-2}$ and
\begin{equation}\label{f0}
F_{\si-1}'(x)=0,\quad x\in J_{\si-1}.
\end{equation}
(The induction process will guarantee the existence of such intervals.) Let $N_{1,\si}$
be such that
\begin{equation}\label{minbj}
\min\{B,|J_{\si-1}|^8\}\ge b_n,\quad n\ge N_{1,\si},
\end{equation}
and let
\begin{equation}\label{n2}
N_{2,\si}:=\left(\frac{\|F^{(2)}_{\si-1}\|}{d_{\si-1}}\right)^{10}.
\end{equation}
Finally, we take
$$
n_{\si}>\max\{n_{\si-1},N_{1,\si},N_{2,\si}\}
$$
so big that the function $f_{n_{\si},b_{n_{\si}}}'$
oscillates a few times inside the interval $J_{\si-1}$
and so it
vanishes on some interval in each oscillation, that is, inside $J_{\si-1}$
there exists an interval
$J_{\si}\subset J_{\si-1}
$
as required in \eqref{f0}.

Now denote
$$
\Phi_{\si}:=\sum_{j=\si}^{\infty}d_{j-1}f_{n_j,b_{n_j}},
$$
where the convergence of the series is justified by the definition of the $d_j$'s and the fact
that $\|f_{n,b_n}\|\le 4\pi,$ for all $n.$ Indeed,
\begin{equation}\label{phisi}
\aligned
\|\Phi_{\si}\|&\le 4\pi d_{\si-1}\left(1+\frac{b^{9/4}_{n_{\si}}}{n_{\si}}+\frac{b^{9/4}_{n_{\si}}b^{9/4}_{n_{\si+1}}}{n_{\si}n_{\si+1}}+...\right) \\
&\le 4\pi d_{\si-1}\sum_{j=0}^{\infty}2^{-j}=8\pi d_{\si-1}.
\endaligned
\end{equation}
Also, since $\|f_{n,b_n}'\|\le 4,$ for all $n,$ then the series is differentiable term by term, which
in turn implies that $\Phi_{\si}\in\Delta^{(1)}(Y^*)$ for all $\si\ge 1.$ So we define
$$
f=f_{\bar\epsilon}:=\sum_{j=1}^{\infty}d_{j-1}f_{n_j,b_{n_j}},
$$
then $f\in\Delta^{(1)}(Y^*)$ and in particular $f$ is monotone in $I_0.$
 Without loss of generality and according with our assumption of even $s$ at the beginning of the section
$f$ is nondecreasing in $I_0$.
{\lem\label{lem1} For each $\si\ge 1$ we have
\begin{equation}\label{om4f}
\omega_4(f,\pi/n_{\si})\le C_8d_{\si}.
\end{equation}
}

\proof
 By \eqref{phisi},
\begin{equation}\label{om4phi}
\omega_4(\Phi_{\si+1},\pi/n_{\si})\le \pi^4 \|\Phi_{\si+1}\|\le 8\pi^5d_{\si}.
\end{equation}
Also, by \eqref{om4},
\begin{equation}\label{om4df}
\omega_4(d_{\si-1}f_{n_{\si},b_{n_{\si}}},\pi/n_{\si})\le d_{\si-1}C_3\frac{b^{9/4}_{n_{\si}}}{n_{\si}}=C_3d_{\si}.
\end{equation}
Finally,
\begin{equation}\label{om4fsi}
\aligned
\omega_4(F_{\si-1},\pi/n_{\si})&\le\pi^2\omega_2(F_{\si-1},\pi/n_{\si})\\
&\le\frac{\pi^4}{n_{\si}^2}\|F^{(2)}_{\si-1}\|\\
&=\pi^4\frac{\|F^{(2)}_{\si-1}\|}{d_{\si-1}}n_{\si}^{-1/10}\left(\frac1{n_{\si}^{2/5}b_{n_{\si}}}\right)^{9/4}d_{\si}\\
&\le \pi^4d_{\si},
\endaligned
\end{equation}
where we used \eqref{n2} and the definitions of $b_{n_{\si}},\ d_{\si}$ and $n_{\si}.$ So \eqref{om4f} follows from combining
\eqref{om4phi}, \eqref{om4df} and \eqref{om4fsi}. Lemma \ref{lem1} is proved.

{\lem\label{lem2} For any measurable $E\subset [-\pi,\pi]$ satisfying
\begin{equation}\label{eep}
|E|\le\epsilon_{n_{\si}},
\end{equation}
and any polynomial $P_{n_{\si}}$ satisfying
\begin{equation}\label{p'nsi}
P_{n_{\si}}'(x)\ge0,\quad x\in I_0\setminus E,
\end{equation}
there exists an absolute constant $C_9$
such that
\begin{equation}\label{f-pnsi}
\|f-P_{n_{\si}}\|\ge (b_{n_{\si}}^{-1/8}/\pi-C_9)d_{\si}.
\end{equation}
}

\proof
Since $F_{\si-1}$ is constant on $J_{\si-1},$
we may write
\begin{equation}\label{fx}
f(x)=d_{\si-1}f_{n_{\si,b_{n_{\si}}}}(x)+\Phi_{\si+1}(x)+M,\quad x\in J_{\si-1}.
\end{equation}
Let
$$
Q_{n_{\si}}:=\frac{1}{d_{\si-1}}(P_{n_{\si}}-M).
$$
Then it follows from \eqref{p'nsi} that
$$
Q_{n_{\si}}'(x)\ge 0,\quad x\in J_{\si-1}\setminus E.
$$
Thus by virtue of Lemma \ref{lem0},
\begin{equation}\label{qn-fn}
\|Q_{n_{\si}}-f_{n_{\si},b_{n_{\si}}}\|_{J_{\si-1}}\ge \frac{b_{n_{\si}}^{2}|J_{\si-1}|}{\pi n_{\si}}-\frac{C_4}{n_{\si}}\left(b_{n_{\si}}^{9/4}+b_{n_{\si}}|E|+
\frac{b_{n_{\si}}^{5/4}}{n_{\si}}\right).
\end{equation}
The definition of $n_{\si}$ and \eqref{minbj} yield
$$
b^2_{n_{\si}}|J_{\si-1}|=b^{17/8}_{n_{\si}}\frac{|J_{\si-1}|}{b^{1/8}_{n_{\si}}}\ge b^{17/8}_{n_{\si}}
$$
On the other hand, \eqref{eep} and the definition of $b_{n_{\si}}$ imply
$$
b_{n_{\si}}|E|\le b_{n_{\si}}\epsilon_{n_{\si}}\le b_{n_{\si}}^{9/4},
$$
and
$$
\frac{b_{n_{\si}}^{5/4}}{n_{\si}}\le b_{n_{\si}}^{15/4}\le b_{n_{\si}}^{9/4}.
$$
Hence \eqref{qn-fn} implies
$$
\|Q_{n_{\si}}-f_{n_{\si},b_{n_{\si}}}\|_{J_{\si-1}}\ge \frac{1}{n_{\si}}\left(b_{n_{\si}}^{17/8}/\pi-3C_4b_{n_{\si}}b_{n_{\si}}^{9/4}\right)=
\frac{b_{n_{\si}}^{9/4}}{n_{\si}}\left(b_{n_{\si}}^{-1/8}/\pi-3C_4\right).
$$
In other words,
$$
\|P_{n_{\si}}-M-d_{\si-1}f_{n_{\si},b_{n_{\si}}}\|_{J_{\si-1}}\ge d_{\si-1}\frac{b_{n_{\si}}^{9/4}}{n_{\si}}\left(b_{n_{\si}}^{-1/8}/\pi-3C_4\right)=
d_{\si}\left(b_{n_{\si}}^{-1/8}/\pi-3C_4\right).
$$
In view of \eqref{fx}, it follows from \eqref{phisi} that,
$$
\aligned
\|f-P_{n_{\si}}\|\ge\|f-P_{n_{\si}}\|_{J_{\si-1}}&\ge \|P_{n_{\si}}-M-d_{\si-1}f_{n_{\si},b_{n_{\si}}}\|_{J_{\si-1}}-\|\Phi_{\si+1}\| \\
&\ge \left(b_{n_{\si}}^{-1/8}/\pi-(3C_4+8\pi)\right)d_{\si},
\endaligned
$$
and Lemma \ref{lem2} is proved with $C_9:=3C_4+8\pi.$

\vskip 0.2cm

\noindent
$\bf 3^\circ .$
The proof of \eqref{main} now follows from Lemmas \ref{lem1} and \ref{lem2}, namely,
$$
\aligned
\limsup_{n\rightarrow\infty}\frac{E_n^{(1)}(f,\epsilon_n,Y^*)}{\omega_4(f,\pi/n)}&\ge \limsup_{\si\rightarrow\infty}\frac{E_{n_{\si}}^{(1)}(f,\epsilon_{n_{\si}},Y^*)}{\omega_4(f,\pi/{n_{\si}})}\\
&\ge \limsup_{\si\rightarrow\infty}\frac{1}{C_8}\left(b_{n_{\si}}^{-1/8}/\pi-C_9\right)=\infty.
\endaligned
$$
Theorem \ref{thm1} is proved.

\vskip 1cm

\footnotesize{

}

\end{document}